\documentclass[12pt]{article}
\usepackage{amssymb}
\newcommand{\Rcc}{\mathbb{R}}
\newcommand{\Ccc}{\mathbb{C}}
\newcommand{\fr}{\frac}
\newcommand{\lb}{\label}

\newcommand{\be}{\begin{equation}}
\newcommand{\ee}{\end{equation}}
\newcommand{\ba}{\begin{array}}
\newcommand{\ea}{\end{array}}
\newcommand{\beqa}{\begin{eqnarray}}

\newcommand{\ka}{\kappa}
\newcommand{\la}{\lambda}

\newcommand{\eeqa}{\end{eqnarray}}

\newcommand{\Ac}{{\cal A}}

\newcommand{\Rc}{{\cal R}}
\newcommand{\Lc}{{\cal L}}

\newcommand{\Ic}{{\cal I}}

\newcommand{\omo}{\omega}
\newcommand{\kd}{\delta}
\newcommand{\CD}{\Delta}
\newcommand{\epl}{\eta_+}
\newcommand{\emi}{\eta_-}
\newcommand{\ot}{\otimes}

\newcommand{\vae}{\varepsilon}
\begin{document}
\title{}
\author{}
\date{}
\begin{flushright}
FGI-99-3 \\

math.QA/9903093
\end{flushright}

\vspace{1cm}

\noindent
{\Large \bf Two Dimensional Fractional Supersymmetry from the
Quantum Poincar\'{e} Group at Roots of Unity}

\vspace{1cm}
\noindent
{\footnotesize Hadji AHMEDOV$^{a,}$}\footnote{E-mail address: 
hagi@gursey.gov.tr}
{\footnotesize and \"{O}mer F. DAYI$^{a,b,}$}\footnote{E-mail 
addresses: dayi@gursey.gov.tr and  dayi@itu.edu.tr.}

\vspace{10pt}

\noindent
$a)${\footnotesize \it Feza G\"{u}rsey Institute,}

\noindent
\hspace{3mm}{\footnotesize \it P.O.Box 6, 81220
\c{C}engelk\"{o}y--Istanbul, Turkey. }

\vspace{10pt}

\noindent
$b)${\footnotesize \it Physics Department, Faculty of Science and
Letters, Istanbul Technical University,}

\noindent
\hspace{3mm}{\footnotesize \it  80626 Maslak--Istanbul,
Turkey.}

\vspace{2cm}
\noindent
{\footnotesize {\bf Abstract.} A group theoretical understanding 
of the two dimensional fractional supersymmetry is given
in terms of the quantum Poincar\'{e} group at roots of unity.
The  fractional supersymmetry algebra and the quantum group
dual to it are presented and the pseudo-unitary, irreducible
representations of them are obtained. The
matrix elements of these representations are explicitly constructed.
}
\newpage

\noindent 
{\bf 1. Introduction}

\vspace{.5cm}

\noindent
Quantum algebras at roots of unity were seen to be
useful in formulations of some physical systems
whose theoretical understanding  is not very clear\cite{om}.
Obviously, to achieve a complete understanding of the
role of these algebras in applications to physical systems one
should know the quantum groups which are dual to them.
Indeed, in \cite{haci}  $E_q(1,1)$ at roots of unity and
in \cite{ha-om} $SL_q(2, \Rcc )$ at roots of unity were constructed.
In the formulation of these groups one is obliged
to introduce some new variables which are the generalized
Grassmannians $\eta_\pm$ satisfying
$\eta_\pm^p=0$ where $p$ is a positive integer.
On the other hand 
these coordinates were
used  to obtain superspace realizations of
the fractional supersymmetry charges\cite{dur}--\cite{dmab}.
Although some algebraic properties of the two dimensional fractional
supersymmetry were discussed  in \cite{ssz},
the correct behaviour under the Lorentz generator
could be obtained
in terms of some  restrictions and a spectral
parameter. 

A group
theoretical understanding of the 
fractional supersymmetry appears to be lacking.
Our aim is to shed some light on the group theoretical aspects  of the
fractional supersymmetry in two dimensions.
Hence, guided by the formulation of the two dimensional
quantum Poincar\'{e} algebra at roots of unity\cite{haci}, 
we introduce the  fractional supersymmetry algebra $U_F$
endowed with a Hopf algebra structure. We present the quantum
group $\Ac_F$ which is its dual. Pseudo-unitary, irreducible 
corepresentations of $\Ac_F$  are given and the matrix elements
of them  are explicitly calculated. We, then, define
the pseudo-unitary quasi-regular corepresentation of $\Ac_F$
and the corresponding $*$-representation of $U_F .$ 

\vspace{1cm}

\noindent
{\bf 2. Fractional supersymmetry algebra and its dual  }

\vspace{.5cm}

\noindent
Let us deal with the two dimensional Poincar\'{e} algebra
$U(e(1,1))$ generated by $P_\pm$ and $H$ possessing the commutators
\be
\lb{f1}
{[ P_+ ,P_- ]} =   0,\ 
{[ P_\pm , H]}   =  \pm i P_\pm ,\\
\ee
and the involutions
\be
\lb{f2}
P_\pm^*=P_\pm  ,\  H^*=H.
\ee
The two dimensional
fractional supersymmetry generators  $p_\pm$ are
defined  to satisfy
\be
\lb{f3}
p_\pm^p=P_\pm ,
\ee
where $p$ is a positive integer, without any condition on
commutation relation of $p_+$ with $p_-$. 
Obviously, the simplest choice is $p_+p_-=p_-p_+.$ 
Thus,  the quantum Poincar\'{e} algebra at roots of
unity $U_q(e(1,1))$ with  $q^p=1,$ $p$ is an odd, positive integer,
generated by $p_\pm$ and $\ka $ satisfying
\be
\lb{f4}
{[ p_+ , p_-]} =0,\ \ka p_\pm   =
q^{\pm 1} p_\pm \ka ,\ \ka^p=1_U,
\ee
with  the involutions
\be
\lb{f5}
p_\pm^*=p_\pm , \ka^*=\ka ,
\ee
suits well with our purposes. $1_U$ denotes the unit element
of the algebra.

The two dimensional fractional supersymmetry algebra
 denoted $U_F \equiv \left(U(e(1,1)),U_q(e(1,1))\right)$ 
is generated by
$P_\pm,\ H,\ p_\pm ,\ \ka$
satisfying  (\ref{f1})--(\ref{f5}) and moreover, the commutation
relations
\be
\lb{fl}
[\ka ,H]=0,\   [p_\pm ,H]=\pm \fr{i}{p} p_\pm .
\ee
The latter is the consequence of (\ref{f1}) and (\ref{f3}).

The basis elements of $U_F$  are 
\be
\phi^{nmkrsl}\equiv p_+^np_-^m\ka^k P_+^rP_-^sH^l ,
\ee
where $n,m,k,r,s,l$ are positive integers.

We can equip  $U_F$ with the Hopf algebra structure 
\be
\begin{array}{lll}
\Delta (P_{\pm })=P_{\pm }\otimes 1_U+1_U\otimes P_{\pm },&
\varepsilon (P_{\pm })=0,&
S(P_{\pm})=-P_{\pm },\\
\CD (H) =H \ot 1_U   +1_U   \ot H, & \vae (H)=0, & S(H)=-H, \\
\Delta (p_\pm )=p_\pm \otimes \kappa +\kappa^{-1}\otimes p_\pm , &
\varepsilon (p_\pm )=0 , &
 S(p_\pm )=-q^{\pm 1}p_\pm , \\
\Delta (\kappa)=\kappa\otimes \kappa, &
\varepsilon (\kappa^{\pm 1})=1, &
S(\kappa^{\pm 1})=\kappa^{\mp 1}.
\end{array}
\ee

Now, we would like to present  the groups which are dual
to the algebras considered above.

The $*$--algebra $\Ac (E(1,1))$ of infinitely differentiable
functions on the two dimensional Poincar\'{e} group
$E(1,1)$ is dual to the algebra $U(e(1,1)).$  For any
$f(z_+,z_-,\la )\in \Ac (E(1,1)) $ we have the involution
\be
z_\pm^*=z_\pm ,\ \la^*=\la , 
\ee
where
\[
\left(
{\begin{array}{ccc}
e^\la & 0 & z_+ \\
0 & e^{-\la} & z_- \\
0 & 0 & 1
\end{array} }
\right)
\in E(1,1) .
\]

The two dimensional quantum Poincar\'{e} group at roots of
unity\cite{haci}
is the $*$--algebra $\Ac (E_q(1,1))$ with $q^p=1,$
generated by $\eta_\pm, \kd $ satisfying
\beqa
\eta _{-}\eta _{+}=q^2\eta _{+}\eta _{-}, & \eta _{\pm }\delta 
=q^2\delta \eta _{\pm }, &  \\
\eta_\pm^p=0 ,& \kd^p =1_A , &  \\
\eta_\pm^*=\eta_\pm,  & \kd^*=\kd  ,&
\eeqa
where $1_A$ is the unit element of $\Ac .$
The dual of $U_F$ is the $*$--algebra
$\Ac_F =\Ac (E(1,1)) \times \Ac (E_q(1,1))$ with the Hopf algebra structure
\[
\begin{array}{lll}
\Delta (\eta_\pm )=\eta_\pm\otimes 1_A+\delta^{\pm 1}
e^{\pm \la /p}\otimes \eta_\pm ,
& \varepsilon (\eta_\pm)=0,
& S(\eta_\pm)=-\delta ^{\mp1} \eta_\pm , \\
\Delta (\delta )=\delta\otimes \delta
& \varepsilon (\delta^{\pm 1} )=1, &
S(\delta^{\pm 1} )=\delta ^{\mp 1}, \\
\CD (\la ) = \la \ot 1_A +1_A \ot \la , &
\vae (\la ) =-\la , &
S (\la )=-\la ,
\end{array}
\]
\[
\Delta (z_{\pm })=z_{\pm }\otimes 1_A+e^{\pm \la }1_A\otimes z_{\pm }
+(-1)^{\frac{p+1}2}\sum_{n=1}^{p-1}
\frac{q^{\pm n^2}}{[p-n]![n]!}\eta _{\pm }^{p-n}
\delta^{\pm n} e^{\pm \la n /p}
\otimes \eta _{\pm }^n ,
\]
\[
S(z_{\pm })=-z_{\pm },\ \ \ \varepsilon (z_{\pm })=0.
\]
We use the symmetric $q$--number
\[
[n] =\fr{q^n-q^{-n}}{q-q^{-1}}
\]
and the $q$--factorial $[n]!=[n][n-1]\cdots [1].$

Since any function of $\Ac (E(1,1))$ can locally be expanded in Taylor series,
there is a local basis of $\Ac_F$ given by
\be
a^{nmktsl}\equiv \eta_+^n \eta_-^m \zeta (k,\kd ) z_+^t z_-^s \la^l ,
\ee
where $n,m,k,t,s,l$ are positive integers and
we defined
$$
\zeta (m, \kd )
\equiv \frac 1p\sum_{n=0}^{p-1}q^{-nm}\delta ^n.
$$
The duality relations between $\Ac_F$ and $U_F$  are
\be
\lb{dr}
\langle \phi^{nmktsl},
a^{n^\prime m^\prime k^\prime t^\prime s^\prime l^\prime } \rangle
 = i^{n+m+t+s+l} q^{\frac{n-m}{2}-nm}l!t!s![n]![m]! 
 \delta_{nn^\prime } \delta_{mm^\prime }
\delta_{tt^\prime }\delta_{ss^\prime } \delta_{ll^\prime }
\delta_{k+n+m, k^\prime }.  
\ee

\vspace{1cm}

\noindent
{\bf 3. Pseudo-unitary, irreducible corepresentations of  $\Ac_F$}

\vspace{.5cm}

\noindent
Let $C_0^\infty (\Rcc )$ be the space of
all infinitely differential functions  of
finite support in $\Rcc$
and $P(t)$ denote the algebra of polynomials in $t$
subject to the conditions $t^p=1$ and $t^*=t.$
The linear map
$$
\pi_r(U_F):  C_0^\infty (\Rcc ) \times P(t)\rightarrow
C_0^\infty (\Rcc ) \times P(t)
$$
given as
\beqa
\pi_r (p_\pm ) f(x) a(t) & = &
(-r)^{1/p} e^{\pm x/p} t^{\pm 1} f(x)a(t) ,\nonumber \\
\pi_r (P_\pm ) f(x) a(t) & = &
-r e^{\pm x} f(x)a(t),  \lb{cc} \\
\pi_r (H ) f(x) a(t) & = &
-i \fr{d}{dx} f(x)a(t), \nonumber \\
\pi_r (\ka ) f(x) a(t) & = &
 f(x)a(qt), \nonumber   
\eeqa
defines the irreducible representation of $U_F$
in $C_0^\infty (\Rcc ) \times P(t) .$

Let us introduce the following Hermitian forms
for the space $C_0^\infty (\Rcc ) \times P(t),$
\beqa
(f_1,f_2) & = & \int_{-\infty}^{+\infty}dxf_1(x)
\overline {f_2(x)} ,\lb{h1} \\
(a_1,a_2) & = & \Phi \left( a_1(t)a^*_2(t) \right), \lb{h2}
\eeqa
where
\be
\Phi(t^s)=\kd_{s,0{\rm (mod\ p)}}.
\ee

$C_0^\infty (\Rcc ) $ endowed with the norm induced by (\ref{h1})
leads to the Hilbert space of the square integrable
functions on $\Rcc .$ On the other hand $P(t)$ with the
pseudo--norm $||a||^2 \equiv (a,a)$ is the pseudo-Euclidean space
with $\fr{p+1}{2}$ positive and $\fr{p-1}{2}$ negative
signatures\cite{haci}. Now, one can verify that 
$\pi_r$ defines pseudo-unitary,
irreducible $*$--representation of $U_F $ for real $r.$

By making use of the duality relations (\ref{dr}),
we can derive from (\ref{cc}) the irreducible
corepresentations of $\Ac_F$ as
\be
\lb{tir}
T_r(f(x)a(t)) = 
\sum_{n,m,k=0}^{p-1}\sum_{t,s,l=0}^{\infty}
\fr{a^{nmktsl} \pi_r(\phi^{nmktsl})f(x)a(t) }
{\langle \phi^{nmktsl},a^{n m k t s l } \rangle},
\ee
which is pseudo-unitary for real $r.$

Consider the Fourier transform of $f(x)\in C_0^\infty (\Rcc ) $
\be
F(\nu )  =  \int_{-\infty}^{+\infty} f(x) e^{\nu x} dx.
\ee
This integral converges for any complex $\nu .$
$F(\nu)$ is an analytic function and moreover, satisfies
$(\nu =\nu_1+i\nu_2)$
\be
|F(\nu_1 +i \nu_2) |< Ke^{c|\nu_1 |},
\ee
for some real constants $K$ and $c.$
Then we can write the  inverse transform as
\be
f(x)  =  \fr{1}{2\pi i} \int_{c-i\infty }^{c+i\infty}
F(\nu )e^{-\nu x} d\nu .
\ee

The Fourier transform of $T_r$ (\ref{tir})
yields the pseudo-unitary corepresentation 
in the space of functions $F(\nu )a(t)$ as
\be
Q_r\left( F(\nu )t^k\right) = 
\int_{c-i\infty }^{c+i\infty} d\mu
\sum_{l=0}^{p-1} Q_{kl}^r(\nu ,\mu , g)
F(\mu )t^l,
\ee
where we denoted the variables as
$g\equiv (g_0;\ g_p)= ( z_+,z_-,\la ;\ \epl ,\emi ,\kd)$
and  the kernel $Q_{kl}^r$ is
\beqa
{\rm for}\ l\geq k, & &   \nonumber \\
Q_{kl}^r(\nu ,\mu ,g )  & = &
(q^{-1/2} \epl )^{l-k} \omo^r_{l-k} (\xi )\kd^k
K_{l-k}^r(\nu ,\mu ,g_0) \nonumber \\
& & + \omo_{p+k-l}^r(\xi )(q^{1/2}\emi )^{p+k-l} \kd^k
K_{l-k-p}^r(\nu ,\mu ,g_0)  \lb{fa1} \\
{\rm for}\ l< k, & &   \nonumber \\
Q_{kl}^r(\nu ,\mu ,g )  & = &
(q^{-1/2} \epl )^{p+l-k} \omo^r_{p+l-k} (\xi )\kd^k
K_{p+l-k}^r(\nu ,\mu ,g_0) \nonumber \\
& & + \omo_{k-l}^r(\xi )(q^{1/2}\emi )^{k-l} \kd^k
K_{l-k}^r(\nu ,\mu ,g_0) .   \lb{fa2} 
\eeqa
We introduced,  in terms of $\xi= q\epl \emi ,$ the polynomials
\be
\omo_s^r (\xi ) =\sum_{m=0}^{p-s-1}
\fr{(ir^{1/p})^{2m+s}}{[n]![m+s]!} (q^s\xi )^m,
\ee
and the functions $K^r_s$ are 
\be
\lb{sf}
K_s^r(\nu ,\mu ,g_0) = \fr{1}{2\pi  i}
e^{\mu \la} \int_{-\infty}^{+\infty}  e^{ ir
( e^xz_+ + e^{-x}z_-)+
x(\nu-\mu+ s/p)} dx.
\ee

By utilizing
the analog of  polar coordinates $\rho >0,\ \beta \in \Rcc $,
the pseudo-Euclidean plane
defined by the axis $z_-=0$ and $z_+=0$ can be studied
in terms of the quadrants
\beqa
{\rm Quad. 1}:&  z_+z_->0 ,\ z_\pm =\fr{1}{2}\rho e^{\pm \beta }, &
{\rm Quad. 2}:\ \ z_+z_-<0 ,\ z_\pm
=\pm \fr{ 1}{2}\rho e^{\pm \beta },  \nonumber   \\
{\rm Quad. 3}:& z_+z_->0,\ z_\pm =\fr{-1}{2}\rho e^{\pm \beta },  &
{\rm Quad. 4}:\ \  z_+z_-<0,\ z_\pm =\mp \fr{ 1}{2}\rho e^{\pm \beta }. 
\nonumber
\eeqa
In these quadrants (\ref{sf}) will lead to the
Hankel functions $H^{(1)}_\nu ,\  H^{(2)}_\nu$ or cylindrical
functions of imaginary argument $K_\nu :$ 
\beqa
{\rm Quad. 1:} &
K_s^r(\nu ,\mu ,g_0)   = 
\fr{1}{2} e^{(\mu -\nu -s/p)( \beta + \fr{\pi i}{2} )
+\mu \la  } H^{(1)}_{\mu-\nu -s/p} (r\rho )  , & \nonumber \\
{\rm Quad. 2:} & K_s^r(\nu ,\mu ,g_0)   = 
\fr{1}{2} e^{(\mu -\nu -s/p)( \beta -\fr{\pi i}{2} )
+\mu \la   }H^{(2)}_{\mu-\nu -s/p} (r\rho ) , & \nonumber  \\
{\rm Quad. 3:} & K_s^r(\nu ,\mu ,g_0)  = 
\fr{1}{\pi i} e^{(\mu -\nu -s/p)( \beta +\fr{\pi i}{2} )
+\mu \la  } K_{\mu-\nu -s/p} (r\rho )  , &\nonumber \\
{\rm Quad. 4:} & K_s^r(\nu ,\mu ,g_0)   = 
\fr{1}{\pi i} e^{ (\mu -\nu -s/p)( \beta -\fr{\pi i}{2} )
+\mu \la  } K_{\mu-\nu -s/p} (r\rho )  , & \nonumber 
\eeqa
with the condition $-1< {\rm Re} (\nu -\mu + s/p)<1.$
\cite{vk}.

\vspace{1cm}

\noindent
{\bf 4. Quasi-regular corepresentation
of $\mathbf \Ac_F$ and $\mathbf *$--representation
of the fractional supersymmetry algebra}

\vspace{.5cm}

\noindent
The comultiplication
\be
\lb{de}
\CD :\Ac \rightarrow \Ac_F \ot \Ac
\ee
defines the pseudo-unitary left quasi-regular corepresentation of
$\Ac_F$ in its subspace $\Ac$ consisting of
the finite sums
$$
X=\sum_s a_s(\epl ,\emi )f_s (z_+,z_-)
$$
where
$a_s(\epl ,\emi )$ are polynomials in $\epl ,\emi$ and
$f_s (z_+,z_-) \in C_0^\infty (\Rcc^2) .$
The space $\Ac$  can be endowed with 
the hermitian form 
\be
\lb{222}
(X,Y)=\Ic_E(XY^*),
\ee
$X,Y\in\Ac$ and
the linear functional
$\Ic_E :\ \Ac \rightarrow \Ccc$
\be
\lb{333}
\Ic_E (X) =\sum_s \Ic (a_s) \Ic_C(f_s)
\ee
is the left invariant integral where\cite{haci}
\beqa
\Ic (\epl^n \emi^m ) & = & q^{-1}\kd_{n, p-1} \kd_{m, p-1}, \\
\Ic_C(f_s)& = & \int_{-\infty}^{+\infty}dz_+dz_-f_s(z_+,z_-).
\eeqa

The right representation of the fractional
supersymmetry algebra $U_F$
corresponding to the quasi-regular
representation (\ref{de}),
\be
\lb{rr}
\Rc (\phi )X= (\phi \ot id)\CD( X),
\ee
$\phi \in U_F,$ is a $*$--representation
\[
(\Rc (\phi )X,Y)_E =(X,\Rc (\phi^* )Y)_E ,
\]
because the hermitian form (\ref{222}) is defined in terms of
the left
invariant integral (\ref{333}).

The right representations
on the variables $\eta_\pm$ and $f(z_+,z_-)$
can explicitly be written as
\be
\lb{reo}
\begin{array}{ll}
{\cal R}( p_\pm ) \eta_\pm^k =i q^{\pm \frac{1}{2}}[k]\eta_\pm^{k-1}, &
{\cal R}( p_\pm ) \eta_\mp^k =0, \\
{\cal R}(\kappa )\eta_\pm^k = q^{\pm k}\eta_\pm^k ,&
\Rc (H)\eta_\pm^n=\pm \fr{in}{p}\eta_\pm^n, \\
{\cal R}( p_\pm ) f=
i q^{\pm \frac{1}{2}}
\frac{(-1)^{\frac{p+1}{2}} }{[p-1]!} \eta^{p-1}_\pm \frac{df}{dz_\pm }, &
{\cal R}( P_\pm )f= i \frac{df}{dz_\pm }, \\
{\cal R}(\kappa ) f= f,  &
\Rc (H) f= iz_+\fr{df}{dz_+}-iz_-\fr{df}{dz_-}.
\end{array}
\ee

In terms of the following relations satisfied by
the right representation ${\cal R}$ 
\beqa
{\cal R}( \phi \phi^\prime ) & = & {\cal R}(\phi^\prime )
 {\cal R}(\phi ) , \nonumber\\
{\cal R}( p_\pm ) (X Y ) & = & {\cal R}( p_\pm )X {\cal R}(\kappa ) Y
+ {\cal R}(\kappa^{-1})X {\cal R}(p_\pm )Y,   \nonumber\\
{\cal R}(\kappa ) (XY) & = & {\cal R}(\kappa )X
{\cal R}(\kappa )Y  \nonumber\\
\Rc (H) (XY) & = & \Rc (H)X Y+ X\Rc (H) Y  \nonumber
\eeqa
we can define the action of an arbitrary
operator ${\cal R}(\phi )$ on any function in $\Ac.$

The quantum algebra which we deal with possesses one
Casimir element $C=p_+p_-.$ As the complete set
of commuting operators we can choose $\Rc (C)$,
$\Rc (H),$  $\Rc (\kappa ) $ and
$\Lc (H),$ $\Lc (\kappa ) $ where
$\Lc (\phi )$ is the left representation of the element $\phi$
defined similar to (\ref{rr}) with the interchange of $\phi$
with the identity $id.$ One can easily observe that
$\Lc (H) X=0$ and  $\Lc (\kappa )X =X$ for any $X\in \Ac ,$
so that, in the space $\Ac$
the matrix elements can
be labeled as $D^r_{n\nu ,m\mu }.$
Indeed, in terms of the kernel $Q^r_{mn}$
 (\ref{fa1}) one observes that
\[
D_{n\nu ,00} = Q_{0n}^r(\nu ,0 ,g ),
\]
$n \in [0,p-1],$ satisfy
\beqa
\Rc (\ka ) D_{n\nu ,00} & = & q^n D_{n\nu ,00}  \nonumber \\
\Rc (H ) D_{n\nu ,00} & = & -i(\nu+n/p) D_{n\nu ,00}  \nonumber \\
\Rc (C ) D_{n\nu ,00} & = & c^2D_{n\nu ,00}  \nonumber \\
\Rc (p_+)D_{n\nu ,00} & = & c
D_{n + 1\nu ,00}, \nonumber \\
\Rc (p_-)D_{n\nu ,00} & = & c
D_{n- 1\nu ,00}, \nonumber \\
\Rc (P_\pm)D_{k\nu ,00} & = & -r
D_{k\nu \pm 1,00} ,\nonumber
\eeqa
where $c=r^{1/p}q$ and we introduced the notation
$ D_{p\nu ,00} \equiv  D_{0\nu +1 ,00} $
and
$ D_{-1\nu ,00} \equiv  D_{p-1\nu -1 ,00} .$

The right representation obtained in (\ref{reo}) can be used to write
the supercharge operators $\Rc (p_\pm )$ in the superspace given by
$\eta_+,\ z_+$ or 
$\eta_-,\ z_-$
 as 
\be
\lb{sps}
\Rc (p_\pm ) =
i q^{\pm \frac{1}{2}} D_\pm^q +
\frac{(-1)^{\frac{p+1}{2}} }{[p-1]!} \eta^{p-1}_\pm \frac{d}{dz_\pm },
\ee
where $D^q_\pm$ are q--derivatives with respect to $\eta_\pm .$
This is the same with the realization of supercharges given 
in \cite{dur}--\cite{dmab}, obtained in terms of $q$--calculus.

\end{document}